\documentclass[12pt]{article}
\usepackage{curves}
\usepackage{amsmath, amssymb, lamsarrow}
\usepackage{amssymb}
\def\mineappendix{
        \setcounter{section}{1}
        \setcounter{subsection}{0}
        \def\thesection{\Alph{section}}
        \def\sectionap{\@startsection  {section}{1}{\z@}
                        {-3.5ex plus-1ex minus-.2ex} {0ex plus.2ex}
                        {\reset@font\Large\bf  Appendix:  \, }
                        }
        }
\makeatother
\def\Proclaim #1. #2\par{\bigbreak\noindent{\sc#1.\enspace}{\it#2}\par}

\newcommand{\gwii}[1]{\left< \hspace{-2pt} \left< \, #1 \,
        \right>  \hspace{-2pt} \right>_{0}}

\newcommand{\gwiione}[1]{\left< \hspace{-2pt} \left< \, #1 \,
        \right> \hspace{-2pt} \right>_{1}}

\newcommand{\gwiitwo}[1]{\left< \hspace{-2pt} \left< \, #1 \,
        \right> \hspace{-2pt} \right>_{2}}

\newcommand{\gwih}[2]{ \left< \, #2 \, \right>_{#1}}
\newcommand{\gwiig}[1]{\left< \hspace{-2pt} \left< \, #1 \,
    \right> \hspace{-2pt} \right>_{g}}
\newcommand{\gwiih}[2]{\left< \hspace{-2pt} \left< \, #2 \,
    \right> \hspace{-2pt} \right>_{#1}}

\newcommand{\grav}[2]{\tau_{#1}(\gamma_{#2})}
\newcommand{\grava}[1]{\tau_{#1}(\gamma_{\alpha})}

\newcommand{\ga}{\gamma_{\alpha}}
\newcommand{\gua}{\gamma^{\alpha}}
\newcommand{\gb}{\gamma_{\beta}}
\newcommand{\gub}{\gamma^{\beta}}
\newcommand{\gm}{\gamma_{\mu}}
\newcommand{\gum}{\gamma^{\mu}}

\newcommand{\vw}{{\cal W}}

\newcommand{\qp}{\circ}

\newtheorem{lem}{Lemma}[section]

\newtheorem{thm}[lem]{Theorem}
\newtheorem{pro}[lem]{Proposition}

\newtheorem{rem}[lem]{Remark}

\newcommand{\M}{\overline{{\cal M}}}
\newcommand{\ba}{{\boldsymbol{\alpha}}}

\newcommand{\cl}{{\cal L}}

\newenvironment{ack}{\noindent \textbf{Acknowledgments}.}

\title{A Genus-3 Topological Recursion Relation}
\author{Takashi Kimura \thanks{Research of the first author was partially supported by NSF grant
DMS-0204824.}, \,\,\, \,\,\, Xiaobo Liu
\thanks{Research of the second author was partially supported by NSF grant
DMS-0505835.}
}
\date{}

\begin{document}
\maketitle

\begin{abstract}
In this paper, we give a new genus-3 topological recursion relation
for Gromov-Witten invariants of compact symplectic manifolds. This
formula also applies to intersection numbers on moduli spaces of
spin curves. A by-product of the proof of this formula is a new
relation in the tautological ring of the moduli space of 1-pointed
genus-3 stable curves.
\end{abstract}

Let $\overline{\cal M}_{g, n}$ be the moduli space of genus-$g$ stable curves with $n$
marked points. It is well known that relations in the tautological rings on
$\overline{\cal M}_{g, n}$ produce universal equations for the Gromov-Witten
invariants of compact symplectic manifolds. Examples of genus-1 and genus-2
universal equations were given in \cite{Ge1}, \cite{Ge2}, and
\cite{BP}. Relations among known universal equations were discussed in
\cite{L2}. It is expected that for manifolds with semisimple quantum
cohomology, such universal equations completely determine all higher genus
Gromov-Witten invariants in terms of its genus-0 invariants. This has been
proven for the genus-1 case in \cite{DZ} and for the genus-2 case in
\cite{L1}. However for genus bigger than 2, no explicit universal equations
had been found except for those which follow from an obvious dimension
count. The main purpose of this paper is to introduce a new genus-3 universal
equation, namely, a genus-3 topological recursion relation.

Associated to the Gromov-Witten invariants of a compact symplectic manifold $M$
is its {\it big phase space}, a product of infinitely many copies of
$H^{*}(M; {\Bbb C})$. We will choose a basis $\{ \ga \mid \alpha = 1, \ldots, N \}$
of $H^{*}(M; {\Bbb C})$. The quantum product $\vw_{1} \qp \vw_{2}$
of two vector fields $\vw_{1}$ and $\vw_{2}$ on the big phase space was introduced
in \cite{L1}. This is an associative product without an identity element. An
operator $T$ on the space of vector fields on the big phase space was also
introduced in \cite{L1} to measure the failure of the string vector field
to be an identity element with respect to this product.
This operator turns out to be a very useful device to translate relations in the
tautological rings of $\overline{\cal M}_{g,n}$ into universal equations for
Gromov-Witten invariants. We will write universal equations of Gromov-Witten
invariants as equations among tensors $\gwiig{\vw_{1} \, \cdots \, \vw_{k}}$
which are defined to be the $k$-th covariant derivatives of the generating
functions of genus-$g$ Gromov-Witten invariants with respect to the trivial
connection on the big phase space. We will briefly review these definitions
in Section \ref{sec:Pre} for completeness.

The main result of this paper is the following theorem.
\begin{thm} \label{thm:g3TRR}
For Gromov-Witten invariants of any compact symplectic manifold,
the following topological recursion relation holds for any vector field $\vw$
on the big phase space:
{\allowdisplaybreaks
\begin{eqnarray}
&& \gwiih{3}{T^{3}(\vw)} \nonumber \\
&=& - \frac{1}{252} \gwiitwo{ \vw \, T(\ga \qp \gua)}
    + \frac{5}{42} \gwiitwo{T(\ga) \, \{ \vw \qp \gua \}} \nonumber \\
&&  + \frac{13}{168} \gwiitwo{ T(\gua) } \gwii{ \ga \, \vw \, \gub \, \gb}
    + \frac{41}{21} \gwiitwo{ T(\gua) } \gwiione{ \left\{\ga \qp \vw\right\}}
        \nonumber \\
&&  - \frac{13}{168} \gwiitwo{ \left\{\vw \qp \ga \qp \gua\right\} }
    + \frac{1}{280} \gwiione{ \vw \gua} \gwiione{ \ga \, \left\{\gub \qp \gb\right\}}
        \nonumber \\
&&  - \frac{23}{5040} \gwiione{\gua} \gwiione{\ga \vw \left\{\gub \qp \gb\right\}}
 - \frac{47}{5040} \gwiione{\gua} \gwiione{\ga \gub} \gwii{\gb \, \vw \, \gum \, \gm}
        \nonumber \\
&&  - \frac{5}{1008} \gwiione{\vw \, \gua} \gwii{\ga \gub \gb \gum} \gwiione{\gm}
  + \frac{23}{504} \gwiione{\gua} \gwii{\ga \vw \gub \gb \gum} \gwiione{\gm}
         \nonumber \\
&&   + \frac{11}{140} \gwiione{\gua \gub} \gwiione{\ga \, \left\{\gb \qp \vw\right\}}
  - \frac{4}{35} \gwiione{\gua} \gwiione{\ga \gub} \gwiione{\left\{\gb \qp \vw\right\}}
            \nonumber \\
&&  + \frac{2}{105} \gwiione{\vw \, \gua} \gwiione{\left\{\ga \qp \gb\right\}} \gwiione{\gub}
  + \frac{89}{210} \gwiione{\gua} \gwii{ \ga \, \vw \, \gub \, \gum} \gwiione{\gb} \gwiione{\gm}
            \nonumber\\
&& - \frac{1}{210} \gwiione{\gua} \gwiione{\ga \gub \left\{\gb \qp \vw\right\}}
    + \frac{1}{140} \gwiione{\vw \, \gua \gub} \gwiione{\left\{\ga \qp \gb\right\}}
            \nonumber \\
&& + \frac{23}{140} \gwiione{\gua \gub} \gwii{\ga \gb \vw \gum} \gwiione{\gm}
    - \frac{3}{140} \gwiione{\gua \gub} \gwiione{\left\{\ga \qp \gb\right\} \vw}
            \nonumber \\
&& - \frac{1}{4480} \gwiione{\vw \gua} \gwii{\ga \gb \gub \gm \gum}
    + \frac{13}{8064} \gwiione{\gua} \gwii{\ga \vw \gub \gb \gum \gm}
    \nonumber \\
&& - \frac{1}{2240} \gwiione{\vw \gua \gub} \gwii{\ga \gb \gum \gm}
    + \frac{41}{6720} \gwiione{\gua \gub} \gwii{\ga \gb \vw \gum \gm}
    \nonumber \\
&& + \frac{1}{53760} \gwii{\vw \gua \ga \gub \gb \gum \gm}
    - \frac{1}{210} \gwiione{\left\{\vw \qp \gua\right\}} \gwiione{\ga \gub \gb}
    \nonumber \\
&& - \frac{1}{5760} \gwiione{\vw \gua \ga \left\{\gub \qp \gb\right\} }
    - \frac{1}{2688} \gwiione{\gua \ga \gub} \gwii{\gb \vw \gum \gm}
    \nonumber \\
&& - \frac{1}{5040} \gwiione{\gua \ga \gub \left\{\gb \qp \vw\right\}}
   + \frac{1}{3780} \gwiione{ \vw \ga \gb \gm } \gwii{\gua \gub \gum}
    \nonumber \\
&&    + \frac{1}{252} \gwiione{ \ga \gb \gm } \gwii{ \vw \gua \gub \gum}.
    \label{eqn:g3TRR}
\end{eqnarray}}
\end{thm}
A by-product of the proof of this theorem is a new relation in
the tautological ring of $\overline{\cal M}_{3, 1}$ which will be given in
Section \ref{sec:taut}. This relation is equivalent to the genus-3
topological recursion relation.

The main idea behind the proof of Theorem~\ref{thm:g3TRR} is that universal
equations for Gromov-Witten invariants can be written as linear combinations
of finitely many terms of a given dimension each corresponding to a boundary
stratum, suitably decorated, in $\overline{\cal M}_{g, n}$.
The goal is to determine the coefficients in these linear combinations.
Suppose that all of the Gromov-Witten invariants of a particular manifold are
known then since these invariants must satisfy the universal relations, it
will imply relations between these coefficients.

In the case of our genus-3 topological recursion relation, the Gromov-Witten
invariants of a point and of ${\mathbb CP}^{1}$, both of which are complete
known, completely determine the coefficients and, hence, this universal
relation.   This method may be adapted  to obtain more universal equations
of any genus and will be explored in a future paper.

Finally, we observe that these universal equations not only apply to Gromov-Witten
theory but also to other cohomological field theories, in the sense of
Kontsevich-Manin, such as r-spin theory \cite{JKV-99} where $r\geq 2$ is an
integer. The correlators in $r$-spin theory are intersection numbers on the
moduli space of $r$-spin curves, (cf. \cite{JKV-99}). $r$-spin theory is
interesting because of the generalized Witten conjecture which states that
its big phase space potential function solves the $r$-th KdV integrable
hierarchy. When $r=2$, this reduces to the ordinary Witten conjecture
\cite{W} proven by Kontsevich \cite{K}. While this conjecture has been proven in
special cases for general $r$ in low genus, it is still open in general.
Since equation \eqref{eqn:g3TRR} also applies to $r$-spin theory, we use it
to calculate some correlators in this theory.

After the completion of this paper, the preprint \cite{AL} appeared in which
a topological recursion relation equivalent to ours on $\M_{3,1}$ was obtained
assuming that the so-called invariance conjectures hold. However, since these
conjectures have not yet been established, their work does not yet prove that
our topological recursion relation holds on $\M_{3,1}$.

\begin{ack}
The authors would like to thank Andreas Gathmann for allowing them to use his
computer program for computing Gromov-Witten invariants of ${\mathbb CP}^{1}$
based on the Virasoro constraints. We would also like to thank the Institut
des Hautes \'Etudes Scientifiques for their hospitality and financial support
where this paper was initiated, and the Focus on Mathematics Program at
Boston University for the usage of their main server whether they are aware
of it or not.
\end{ack}

\section{Preliminaries}
\label{sec:Pre}

Let $M$ be a compact symplectic manifold.
The {\it big phase space} is by definition the infinite product
\[  P := \prod_{n=0}^{\infty} H^{*}(M; {\Bbb C}). \]
Fix a basis $\{ \gamma_{0}, \ldots, \gamma_{N} \}$ of
$H^{*}(M; {\Bbb C})$, where $\gamma_{0}$ is the identity element, of the ordinary
cohomology ring of $M$. Then we denote the corresponding basis for
the $n$-th copy of $H^{*}(M; {\Bbb C})$ in $P$ by
$\{\tau_{n}(\gamma_{0}), \ldots, \tau_{n}(\gamma_{N}) \}$.
We call $\grava{n}$ a {\it descendant of $\gamma_{\alpha}$ with descendant
level $n$}.
We can think of $P$ as an infinite dimensional vector space with a basis
$\{ \grava{n} \mid 0 \leq \alpha \leq N, \, \, \, n \in {\Bbb Z}_{\geq 0} \}$
where ${\Bbb Z}_{\geq 0} = \{ n \in {\Bbb Z} \mid n \geq 0\}$.
Let
$(t_{n}^{\alpha} \mid 0 \leq \alpha \leq N, \, \, \, n \in {\Bbb Z}_{\geq 0})$
be the corresponding coordinate system on $P$.
For convenience, we identify $\grava{n}$ with the coordinate vector field
$\frac{\partial}{\partial t_{n}^{\alpha}}$ on $P$ for $n \geq 0$.
If $n<0$, $\grava{n}$ is understood to be the $0$ vector field.
We also abbreviate $\grava{0}$ by $\gamma_{\alpha}$.
We use $\tau_{+}$ and $\tau_{-}$ to denote the operators which shift the level
of descendants by $1$, i.e.
\[ \tau_{\pm} \left(\sum_{n, \alpha} f_{n, \alpha} \grava{n}\right)
    = \sum_{n, \alpha} f_{n, \alpha} \grava{n \pm 1} \]
where $f_{n, \alpha}$ are functions on the big phase space.

We will adopt the following {\it notational conventions}:
Lower case Greek letters, e.g. $\alpha$, $\beta$, $\mu$, $\nu$,
$\sigma$,..., etc., will be used to index the cohomology classes on $M$.
These indices run from $0$ to $N$. Lower case English
letters, e.g. $i$, $j$, $k$, $m$, $n$, ..., etc., will be used to
index the level of descendants. These indices run over the set of all
non-negative integers, i.e. ${\Bbb Z}_{\geq 0}$. All summations are
over the entire ranges of the corresponding indices unless otherwise
indicated.
Let
\[ \eta_{\alpha \beta} = \int_{M} \gamma_{\alpha} \cup
    \gamma_{\beta}
\]
 be the intersection form on $H^{*}(M, {\Bbb C})$.
We will use $\eta = (\eta_{\alpha \beta})$ and $\eta^{-1} =
(\eta^{\alpha \beta})$ to lower and raise indices.
For example,
\[ \gua :=  \eta^{\alpha \beta} \gb.\]
Here we are using the summation convention that repeated
indices (in this formula, $\beta$) should be summed
over their entire ranges.

Let
\[ \gwih{g, d}{\grav{n_{1}}{\alpha_{1}} \, \grav{n_{2}}{\alpha_{2}} \,
    \ldots \, \grav{n_{k}}{\alpha_{k}}} :=
    \int_{[\overline{\cal M}_{g,n}(M;d)]^\mathrm{virt}} \bigcup_{i=1}^k
    (\Psi^{n_{i}} \cup \mathrm{ev_{i}}^*\gamma_{\alpha_{i}})
\]
be the genus-$g$, degree $d$,
 descendant Gromov-Witten invariant associated
to $\gamma_{\alpha_{1}}, \ldots, \gamma_{\alpha_{k}}$ and nonnegative
integers $n_{1}, \ldots, n_{k}$
(cf. \cite{W}, \cite{RT}, \cite{LiT}).
Here, $\overline{\cal M}_{g, k}(M; d)$ is the moduli space of stable maps
from genus-$g$, $k$-pointed curves to $M$ of degree $d \in H_{2}(M; {\Bbb Z})$.
$\Psi_{i}$ is the first Chern class of the tautological line bundle
over $\overline{\cal M}_{g, k}(M; d)$ whose geometric fiber over a stable map
is the cotangent space of the domain curve at the $i$-th marked point while
$\mathrm{ev_i}:\overline{\cal M}_{g,n}(M; d)\to M$ is the $i$-th evaluation map for all
 $i=1,\ldots,k$. Finally, $[\overline{\cal M}_{g,n}(M;d)]^\mathrm{virt}$ is the virtual
 fundamental class.
The genus-$g$
generating function is defined to be
\[ F_{g} =  \sum_{k \geq 0} \frac{1}{k!}
         \sum_{ \begin{array}{c}
        {\scriptstyle \alpha_{1}, \ldots, \alpha_{k}} \\
                {\scriptstyle  n_{1}, \ldots, n_{k}}
                \end{array}}
                t^{\alpha_{1}}_{n_{1}} \cdots t^{\alpha_{k}}_{n_{k}}
    \sum_{d \in H_{2}(V, {\Bbb Z})} q^{d}
    \gwih{g, d}{\grav{n_{1}}{\alpha_{1}} \, \grav{n_{2}}{\alpha_{2}} \,
        \ldots \, \grav{n_{k}}{\alpha_{k}}} \]
where $q^{d}$ belongs to the Novikov ring.
This function is understood as a formal power series n the variables
    $\{\,t_{n}^{\alpha}\,\}$ with coefficients in the Novikov ring.

Introduce a $k$-tensor
 $\left< \left< \right. \right. \underbrace{\cdot \cdots \cdot}_{k}
        \left. \left. \right> \right> $
defined by
\[ \gwiig{{\cal W}_{1} {\cal W}_{2} \cdots {\cal W}_{k}} \, \,
         := \sum_{m_{1}, \alpha_{1}, \ldots, m_{k}, \alpha_{k}}
                f^{1}_{m_{1}, \alpha_{1}} \cdots f^{k}_{m_{k}, \alpha_{k}}
        \, \, \, \frac{\partial^{k}}{\partial t^{\alpha_{1}}_{m_{1}}
            \partial t^{\alpha_{2}}_{m_{k}} \cdots
            \partial t^{\alpha_{k}}_{m_{k}}} F_{g},
 \]
for vector fields ${\cal W}_{i} = \sum_{m, \alpha}
        f^{i}_{m, \alpha} \, \frac{\partial}{\partial t_{m}^{\alpha}}$ where
$f^{i}_{m, \alpha}$ are functions on the big phase space.
This tensor is called the {\it $k$-point (correlation) function}.

For any vector fields $\vw_{1}$ and $\vw_{2}$ on the big phase space,
the quantum product of $\vw_{1}$ and $\vw_{2}$ is defined by
\[ \vw_{1} \qp \vw_{2} := \gwii{\vw_{1} \, \vw_{2} \, \gua} \ga. \]
Define the vector field
\[ T(\vw) := \tau_{+}(\vw) - \gwii{\vw \, \gua} \ga \]
for any vector field $\vw$. The operator $T$
was introduced in \cite{L1} as a convenient tool in the study of universal
equations for Gromov-Witten invariants. Let $\psi_{i}$ be the first
Chern class of the tautological line bundle over $\overline{\cal M}_{g, k}$
whose geometric fiber over a stable curve is the cotangent space of the curve
at the $i$-th marked point.  When we translate a relation in the tautological
ring of $\overline{\cal M}_{g, k}$ into differential equations for
generating functions of Gromov-Witten invariants,
each $\psi$ class corresponds to the insertion of the operator $T$. Let
$\nabla$ be the trivial flat connection on the big phase space with respect
to which $\grava{n}$ are parallel vector fields for all $\alpha$ and
$n$. Then the covariant derivative of the quantum product satisfies
\[ 
\nabla_{\vw_{3}} (\vw_{1} \qp \vw_{2})
    = (\nabla_{\vw_{3}} \vw_{1}) \qp \vw_{2}
        +  \vw_{1} \qp (\nabla_{\vw_{3}} \vw_{2})
        + \gwii{\vw_{1} \, \vw_{2} \, \vw_{3} \, \gua} \ga
\]
and the covariant derivative of the operator $T$ is given by
\[ 
 \nabla_{\vw_{2}} \, \, T(\vw_{1}) = T(\nabla_{\vw_{2}} \vw_{1})
    - \vw_{2} \qp \vw_{1}
\]
for any vector fields $\vw_{1}, \vw_{2}$ and  $\vw_{3}$
(cf. \cite[Equation (8) and Lemma 1.5]{L1}).
We need to use these formulas in order to compute derivatives of universal equations.

\section{Proof of the genus-3 topological recursion relation}

The cohomology class $\psi_{1}^{g}$ vanishes on ${\cal M}_{g,1}$ due
to a result of Ionel (cf. \cite{Io}). It was proven in \cite{GV} and
\cite{Io} that $\psi_{1}^{g}$ on $\overline{\cal M}_{g,1}$ is
supported on the locus of curves which has at least one genus-0
component. Furthermore, by a result of Faber and Pandharipande
\cite{FP}, $\psi_{1}^{g}$ is equal to a  class from the boundary
strata which is tautological, and therefore is a linear combination
of products of $\psi$ and $\kappa$ classes and fundamental classes
of some boundary strata. For $g=3$, $\kappa$ classes do not occur in
this linear combination since components of curves in the boundary
strata have genus at most 2 and $\kappa_{1}$ can be represented as
linear combinations of $\psi$ classes and fundamental classes of
boundary strata  on the moduli spaces of stable curves of genus less
than or equal to 2 (cf. \cite{AC}). Therefore, it follows that
$\psi_{1}^{3}$ on $\overline{\cal M}_{3,1}$ can be written as a
linear combination of products of the $\psi$ classes and the
fundamental classes of some boundary strata. By taking into
consideration the genus-0 and genus-1 topological recursion
relations as well as Mumford's genus-2  relation, we can translate
these results into the following universal equations with unknown
constants $a_{1}, \ldots, a_{30}$: {\allowdisplaybreaks
\begin{eqnarray}
0&=& \Phi(\vw) \nonumber \\
&:=& - \gwiih{3}{T^{3}(\vw)} \nonumber \\
&& +  a_{1} \gwiitwo{ T(\vw) \left\{\ga \qp \gua\right\}}
    + a_{2} \gwiitwo{ \vw \, T(\ga \qp \gua)} \nonumber \\
&&  + a_{3} \gwiitwo{ T(\gua) } \gwii{ \ga \, \vw \, \gub \, \gb}
    + a_{4} \gwiitwo{ T(\gua) } \gwiione{ \left\{\ga \qp \vw\right\}} \nonumber \\
&&  + a_{5} \gwiitwo{ \left\{\vw \qp \ga \qp \gua\right\} }
    + a_{6} \gwiione{ \vw \gua} \gwiione{ \ga \, \left\{\gub \qp \gb\right\}} \nonumber\\
&&  + a_{7} \gwiione{\gua} \gwiione{\ga \vw \left\{\gub \qp \gb\right\}}
    + a_{8} \gwiione{\gua} \gwiione{\ga \gub} \gwii{\gb \, \vw \, \gum \, \gm} \nonumber \\
&&  + a_{9} \gwiione{\vw \, \gua} \gwii{\ga \gub \gb \gum} \gwiione{\gm} \nonumber \\
&&  + a_{10} \gwiione{\gua} \gwii{\ga \vw \gub \gb \gum} \gwiione{\gm}
    + a_{11} \gwiione{\gua \gub} \gwiione{\ga \, \left\{\gb \qp \vw\right\}} \nonumber \\
&& + a_{12} \gwiione{\gua} \gwiione{\ga \gub} \gwiione{\left\{\gb \qp \vw\right\}} \nonumber \\
&&  + a_{13} \gwiione{\vw \, \gua} \gwiione{\left\{\ga \qp \gb\right\}} \gwiione{\gub} \nonumber \\
&& + a_{14} \gwiione{\gua} \gwii{ \ga \, \vw \, \gub \, \gum} \gwiione{\gb} \gwiione{\gm} \nonumber \\
&& + a_{15} \gwiione{\gua} \gwiione{\ga \gub \left\{\gb \qp \vw\right\}}
    + a_{16} \gwiione{\vw \, \gua \gub} \gwiione{\left\{\ga \qp \gb\right\}} \nonumber \\
&& + a_{17} \gwiione{\gua \gub} \gwii{\ga \gb \vw \gum} \gwiione{\gm}
    + a_{18} \gwiione{\gua \gub} \gwiione{\left\{\ga \qp \gb\right\} \vw} \nonumber \\
&& + a_{19} \gwiione{\vw \gua} \gwii{\ga \gb \gub \gm \gum}
    + a_{20} \gwiione{\gua} \gwii{\ga \vw \gub \gb \gum \gm} \nonumber \\
&& + a_{21} \gwiione{\vw \gua \gub} \gwii{\ga \gb \gum \gm}
    + a_{22} \gwiione{\gua \gub} \gwii{\ga \gb \vw \gum \gm} \nonumber \\
&& + a_{23} \gwii{\vw \gua \ga \gub \gb \gum \gm}
    + a_{24} \gwiione{\left\{\vw \qp \gua\right\}} \gwiione{\ga \gub \gb} \nonumber \\
&& + a_{25} \gwiione{\vw \gua \ga \left\{\gub \qp \gb\right\} }
    + a_{26} \gwiione{\gua \ga \gub} \gwii{\gb \vw \gum \gm} \nonumber \\
&& + a_{27} \gwiione{\gua \ga \gub \left\{\gb \qp \vw\right\}}
    + a_{28} \gwiitwo{T(\ga) \, \{ \vw \qp \gua \}} \nonumber \\
&&    + a_{29} \gwiione{ \vw \ga \gb \gm } \gwii{\gua \gub \gum}
    + a_{30} \gwiione{ \ga \gb \gm } \gwii{ \vw \gua \gub \gum} \label{eqn:g3TRRa}
\end{eqnarray}}
where $\vw$ is any vector field on the big phase space.

Using the genus-2 equation discovered by Belurousski-Pandharipande \cite{BP}, we can
write
\[ \gwiitwo{ T(\vw) \left\{\ga \qp \gua\right\}} \]
  as a linear combination
of other terms on the right hand side of equation \eqref{eqn:g3TRRa}
(cf. \cite[Equation (21)]{L1}).
Therefore,  we can set
\begin{equation}  a_{1} = 0. \label{eqn:a1}
\end{equation}
Note that equation \eqref{eqn:g3TRRa} holds for any compact symplectic manifold.
However, we shall see that the Gromov-Witten invariants of a point and
of ${\mathbb CP}^{1}$ already completely determine the coefficients $a_{2},
\ldots, a_{30}$.

\subsection{Relations obtained from the Gromov-Witten invariants of a point}
When the target manifold is a point, all stable maps must necessarily
have degree $0$. Hence, we will omit any reference to the degrees of these
Gromov-Witten invariants. In fact, the moduli space of stable maps into a point
is isomorphic to the moduli space of stable curves.

Since the cohomology space of a point is one dimensional,
coordinates on the big phase space are simply denoted by $t_{0},
t_{1}, t_{2}, \cdots$. We also identify vector fields $
\frac{\partial}{\partial t_{m}}$ with $\tau_{m}$ on the big phase
space.

Gromov-Witten invariants of a point obey the string equation
\begin{eqnarray*}
\gwih{g}{\tau_{0} \, \tau_{n_{1}} \, \cdots \, \tau_{n_{k}}}
    &=& \sum_{j=1}^{k}
        \gwih{g, d}{\tau_{n_{1}} \, \cdots \, \tau_{n_{j}-1}
            \, \cdots \, \tau_{n_{k}}}
        + \delta_{g, 0} \delta_{k, 2} \delta_{n_{1}, 0} \delta_{n_{2}, 0}
\end{eqnarray*}
and the dilaton equation
\begin{eqnarray*}
\gwih{g}{\tau_{1} \, \tau_{n_{1}} \, \cdots \, \tau_{n_{k}}}
    &=& (2g-2+k)
        \gwih{g}{\tau_{n_{1}} \, \cdots  \, \tau_{n_{k}}}
                + \frac{1}{24} \, \delta_{g, 1} \delta_{k, 0}.
\end{eqnarray*}
We can use these two equations to compute Gromov-Witten invariants involving
only $\tau_{0}$ and $\tau_{1}$. More complicated Gromov-Witten invariants
of a point can be computed using the Virasoro constraints, or equivalently
the KdV hierarchy which were conjectured by Witten and proven by
Kontsevich (cf. \cite{W} and \cite{K}). To determine the
coefficients $a_{2}, \ldots, a_{30}$ in equation \eqref{eqn:g3TRR}, we
will only need the genus-2 invariants
\[
\gwih{2}{\tau_{4}} = \frac{1}{1152},
              \hspace{60pt}
\gwih{2}{\tau_{2} \tau_{3}} = \frac{29}{5760},
               \hspace{60pt}
\gwih{2}{\tau_{2} \tau_{2} \tau_{2}} = \frac{7}{240}
\]
and the genus-3 invariants
\[
\gwih{3}{\tau_{7}} =  \frac{1}{82944},
               \hspace{60pt}
\gwih{3}{\tau_{2} \tau_{6} } = \frac{77}{414720},
        \hspace{60pt}
\gwih{3}{\tau_{3} \tau_{5} } = \frac{503}{1451520},
 \]
\[
\gwih{3}{\tau_{4} \tau_{4} } =  \frac{607}{1451520},
             \hspace{40pt}
\gwih{3}{\tau_{2} \tau_{2} \tau_{5} } = \frac{17}{5760},
        \hspace{60pt}
\gwih{3}{\tau_{2} \tau_{3} \tau_{4} } = \frac{1121}{241920},
\]
\[
\gwih{3}{\tau_{3} \tau_{3} \tau_{3} } = \frac{583}{96768},
           \hspace{40pt}
\gwih{3}{\tau_{2} \tau_{2} \tau_{2} \tau_{4} } = \frac{53}{1152},
             \hspace{50pt}
\gwih{3}{\tau_{2} \tau_{2} \tau_{3} \tau_{3} } = \frac{205}{3456},
          \]
\[
\gwih{3}{\tau_{2} \tau_{2} \tau_{2} \tau_{2} \tau_{3} } = \frac{193}{288}.
\hspace{300pt}
\]
All other invariants needed to determine the coefficients $a_{2}, \ldots,
a_{30}$ can be computed from the string and dilaton equations.

We will compute
derivatives of $\Phi(\tau_{m})$ restricted to the origin $t=0$ of the big
phase space. By equation \eqref{eqn:g3TRRa}, these values must all be equal to zero.
We thus obtain some linear relations between the coefficients $a_{2}, \ldots,
a_{30}$ in equation \eqref{eqn:g3TRR}.

From
 $\Phi(\tau_{4})\mid_{t=0} = 0$, we obtain
\begin{equation} 0 = -\frac{1}{82944} + \frac{a_{2}}{384} + a_{23}
    + \frac{a_{25}}{24} + \frac{a_{29}}{24}.
    \label{eqn:R1}
    \end{equation}
From $\tau_{5} \Phi(\tau_{0}) \mid _{t=0} = 0$, we obtain
\begin{equation} 0= -\frac{503}{1451520} + \frac{a_{2}}{288} + \frac{a_{5}}{1152} + a_{23} +
      \frac{a_{25}}{24} + \frac{a_{27}}{24} + \frac{a_{28}}{288} + \frac{a_{29}}{24}.
      \end{equation}
From $\tau_{4} \Phi(\tau_{1}) \mid _{t=0} = 0$, we obtain
\begin{equation} 0= -\frac{607}{1451520} + \frac{a_{2}}{96} + \frac{a_{3}}{384} + 5 a_{23} +
      \frac{a_{25}}{6} + \frac{a_{26}}{24} + \frac{a_{29}}{6} + \frac{a_{30}}{24}.
      \end{equation}
From $\tau_{3} \Phi(\tau_{2}) \mid _{t=0} = 0$, we obtain
\begin{equation}0 = -\frac{503}{1451520} + \frac{29 a_{2}}{1440} + \frac{a_{6}}{576} +
      \frac{a_{18}}{576} + \frac{a_{19}}{24} + \frac{a_{22}}{24} +
      10 a_{23} + \frac{7 a_{25}}{24} + \frac{7 a_{29}}{24}.
      \end{equation}
From $\tau_{2} \Phi(\tau_{3}) \mid _{t=0} = 0$, we obtain
\begin{equation}0 = -\frac{77}{414720} + \frac{29 a_{2}}{1440} + \frac{a_{7}}{576} + \frac{a_{16}}{576} +
      \frac{a_{20}}{24} + \frac{a_{21}}{24} +
      10 a_{23} + \frac{7 a_{25}}{24} + \frac{7 a_{29}}{24}.
      \end{equation}
From $\tau_{3} \tau_{3} \Phi(\tau_{0}) \mid _{t=0} = 0$, we obtain
\begin{eqnarray}0 &=& -\frac{583}{96768} + \frac{29 a_{2}}{576} + \frac{29 a_{5}}{2880} +
      \frac{a_{6}}{288} + \frac{a_{11}}{288} + \frac{a_{18}}{288} + \frac{a_{19}}{12}
      + \frac{a_{22}}{12} \nonumber \\
      && +
      20 a_{23} + \frac{7 a_{25}}{12} + \frac{7 a_{27}}{12} + \frac{29 a_{28}}{576}
      + \frac{7 a_{29}}{12}.
      \end{eqnarray}
From $\tau_{2} \tau_{4} \Phi(\tau_{0}) \mid _{t=0} = 0$, we obtain
\begin{eqnarray}0 &=& -\frac{1121}{241920} + \frac{11 a_{2}}{288} + \frac{a_{3}}{384} +
      \frac{a_{4}}{9216} + \frac{11 a_{5}}{1440} + \frac{a_{7}}{576} + \frac{a_{15}}{576}
      + \frac{a_{16}}{576} +
      \frac{a_{20}}{24} \nonumber \\
      && + \frac{a_{21}}{24} + 15 a_{23} + \frac{a_{24}}{576} + \frac{11 a_{25}}{24} +
      \frac{a_{26}}{24} + \frac{11 a_{27}}{24} + \frac{11 a_{28}}{288} + \frac{11 a_{29}}{24}
        + \frac{a_{30}}{24}.
        \end{eqnarray}
From $\tau_{2} \tau_{3} \Phi(\tau_{1}) \mid _{t=0} = 0$, we obtain
\begin{eqnarray}0  &=& - \frac{1121}{241920} + \frac{29 a_{2}}{288} + \frac{29 a_{3}}{1440}
    + \frac{a_{6}}{288} + \frac{a_{7}}{192} + \frac{a_{8}}{576} + \frac{a_{16}}{192}
    + \frac{a_{17}}{576} + \frac{a_{18}}{288} \nonumber\\
    && + \frac{a_{19}}{12} + \frac{a_{20}}{6} + \frac{a_{21}}{8} + \frac{a_{22}}{8} +
      60 a_{23} + \frac{35 a_{25}}{24} + \frac{7 a_{26}}{24} + \frac{35 a_{29}}{24}
      + \frac{7 a_{30}}{24}.
      \end{eqnarray}
From $\tau_{2} \tau_{2} \Phi(\tau_{2}) \mid _{t=0} = 0$, we obtain
\begin{eqnarray}0 &=& -\frac{17}{5760} + \frac{7 a_{2}}{48} + \frac{a_{6}}{144} + \frac{a_{7}}{72} +
      \frac{a_{9}}{288} + \frac{a_{10}}{288} + \frac{a_{13}}{6912} + \frac{a_{16}}{72} \nonumber \\
      && + \frac{a_{18}}{144} +
      \frac{a_{19}}{4} + \frac{a_{20}}{2} + \frac{a_{21}}{3} + \frac{a_{22}}{6} + 90 a_{23} + 2 a_{25} +
      2 a_{29}.
      \end{eqnarray}
From $\tau_{2} \tau_{2} \tau_{3} \Phi(\tau_{0}) \mid _{t=0} = 0$, we obtain
\begin{eqnarray}0 &=&  -\frac{205}{3456} + \frac{5 a_{2}}{12} + \frac{29 a_{3}}{720}
    + \frac{29 a_{4}}{17280} + \frac{5 a_{5}}{72} + \frac{a_{6}}{72} + \frac{7 a_{7}}{288} +
      \frac{a_{8}}{288} + \frac{a_{9}}{288} + \frac{a_{10}}{288} \nonumber \\
      && + \frac{a_{11}}{72}
      + \frac{a_{12}}{6912} +
      \frac{a_{13}}{6912} + \frac{7 a_{15}}{288} + \frac{7 a_{16}}{288} + \frac{a_{17}}{288} +
      \frac{a_{18}}{72} + \frac{5 a_{19}}{12} + \frac{5 a_{20}}{6}  + \frac{7 a_{21}}{12}
      + \frac{ 5 a_{22}}{12} \nonumber \\
      && + 210 a_{23}  + \frac{7 a_{24}}{288} + \frac{59 a_{25}}{12} + \frac{7 a_{26}}{12} +
    \frac{59 a_{27}}{12} + \frac{5 a_{28}}{12} + \frac{59 a_{29}}{12} + \frac{7 a_{30}}{12}.
    \end{eqnarray}
From $\tau_{2} \tau_{2} \tau_{2} \Phi(\tau_{1}) \mid _{t=0} = 0$, we obtain
\begin{eqnarray}0 &=& -\frac{53}{1152} + \frac{7 a_{2}}{8} + \frac{7 a_{3}}{48}
    + \frac{a_{6}}{24} +
      \frac{a_{7}}{12} + \frac{a_{8}}{48} + \frac{a_{9}}{48} + \frac{a_{10}}{32} + \frac{a_{13}}{1152}
       + \frac{a_{14}}{2304} + \frac{a_{16}}{12} \nonumber \\
      && + \frac{a_{17}}{48} +
      \frac{a_{18}}{24} + \frac{3 a_{19}}{2} + \frac{15 a_{20}}{4} +
      2 a_{21} + \frac{3 a_{22}}{2} + 630 a_{23} \nonumber \\
      && + 12 a_{25} + 2 a_{26} +
      12 a_{29} + 2 a_{30}.
      \end{eqnarray}
From $\tau_{2} \tau_{2} \tau_{2} \tau_{2} \Phi(\tau_{0}) \mid _{t=0} = 0$, we obtain
\begin{eqnarray} 0 &=& -\frac{193}{288} + \frac{49 a_{2}}{12} + \frac{7 a_{3}}{12}
    + \frac{7 a_{4}}{288} + \frac{7 a_{5}}{12} + \frac{a_{6}}{6} + \frac{a_{7}}{3} + \frac{a_{8}}{12}
     + \frac{a_{9}}{12} +
      \frac{a_{10}}{8} + \frac{a_{11}}{6} \nonumber \\
      && + \frac{a_{12}}{288} + \frac{a_{13}}{288}
      + \frac{a_{14}}{576} + \frac{a_{15}}{3} +
      \frac{a_{16}}{3} + \frac{a_{17}}{12} + \frac{a_{18}}{6} + 6 a_{19} + 15 a_{20} + 8 a_{21} +
      6 a_{22} \nonumber \\
      && + 2520 a_{23} + \frac{a_{24}}{3} + 48 a_{25} + 8 a_{26} +
      48 a_{27} + \frac{49 a_{28}}{12} + 48 a_{29} + 8 a_{30}.
      \end{eqnarray}

\subsection{Relations obtained from the Gromov-Witten invariants of ${\mathbb CP}^{1}$}

When the target manifold is ${\mathbb CP}^{1}$, the degrees of the stable
  maps are indexed by $H_{2}({\mathbb CP}^{1}; {\mathbb Z}) \cong {\mathbb
  Z}$. The degree $d$ part of any equation for generating functions of
  Gromov-Witten invariants is the coefficient of $q^{d}$ in the Novikov ring.
We choose the basis $\{ \gamma_{0}, \gamma_{1} \}$ for $H^{*}({\mathbb CP}^{1}; {\mathbb C})$ with
$\gamma_{0} \in H^{0}({\mathbb CP}^{1}; {\mathbb C})$ being the identity of the ordinary
cohomology ring and $\gamma_{1} \in H^{2}({\mathbb CP}^{1}; {\mathbb C})$ the Poincare dual
to a point. Coordinates on the big phase space are denoted by
$\{ t_{n}^{0}, t_{n}^{1} \mid n \in {\mathbb Z}_{+} \}$. We identify
vector fields $\frac{\partial}{\partial t_{n}^{0}}$ and $\frac{\partial}{\partial t_{n}^{1}}$
with $\tau_{n, 0}$ and $\tau_{n, 1}$ respectively. We also define $\tau_{n, \alpha} =0$
if $n < 0$.

The Gromov-Witten invariants of ${\mathbb CP}^{1}$ obey three basic equations:
the string equation
\begin{eqnarray*}
\gwih{g, d}{\tau_{0,0} \, \tau_{n_{1}, \alpha_{1}} \, \cdots \, \tau_{n_{k}, \alpha_{k}}}
    &=& \sum_{j=1}^{k}
        \gwih{g, d}{\tau_{n_{1}, \alpha_{1}} \, \cdots \, \tau_{n_{j}-1, \alpha_{j}}
            \, \cdots \, \tau_{n_{k}, \alpha_{k}}}  \\
    &&        + \delta_{g, 0} \delta_{d, 0} \delta_{k, 2} \delta_{n_{1}, 0} \delta_{n_{2}, 0}
            \delta_{\alpha_{1} + \alpha_{2}, 1},
\end{eqnarray*}
the dilaton equation
\begin{eqnarray*}
\gwih{g, d}{\tau_{1,0} \, \tau_{n_{1}, \alpha_{1}} \, \cdots \, \tau_{n_{k}, \alpha_{k}}}
    &=& (2g-2+k)
        \gwih{g, d}{\tau_{n_{1}, \alpha_{1}} \, \cdots  \, \tau_{n_{k}, \alpha_{k}}}
                + \frac{1}{12} \, \delta_{g, 1} \delta_{k, 0} \delta_{d, 0},
\end{eqnarray*}
and the divisor equation
\begin{eqnarray*}
\gwih{g, d}{\tau_{0,1} \, \tau_{n_{1}, \alpha_{1}} \, \cdots \, \tau_{n_{k}, \alpha_{k}}}
    &=& d \gwih{g, d}{\tau_{n_{1}, \alpha_{1}} \,\cdots \,  \tau_{n_{k}, \alpha_{k}}}  \\
    && + \sum_{j=1}^{k} \delta_{\alpha_{j}, 0}
            \gwih{g, d}{\tau_{n_{1}, \alpha_{1}} \, \cdots \, \tau_{n_{j}-1, 1}
            \, \cdots \, \tau_{n_{k}, \alpha_{k}}}  \\
    &&        + \delta_{g, 0} \delta_{d, 0} \delta_{k, 2} \delta_{n_{1}, 0} \delta_{n_{2}, 0}
            \delta_{\alpha_{1}, 0} \delta_{\alpha_{2}, 0}
            - \frac{1}{24} \, \delta_{g, 1} \delta_{d, 0} \delta_{k, 0} .
\end{eqnarray*}
We can use these three equations to compute Gromov-Witten invariants
involving only $\tau_{0,0}$, $\tau_{1,0}$, and $\tau_{0,1}$. More complicated
Gromov-Witten invariants for ${\mathbb CP}^{1}$  can be computed using
the Virasoro constraints which was conjectured in
\cite{EHX} and proven in \cite{Gi}.
A computer program for computing such invariants based on
the Virasoro constraints was written by Andreas Gathmann
(cf. \cite{Ga}). To determine $a_{2}, \cdots a_{30}$ in equation
\eqref{eqn:g3TRR}, we only need a small number of such invariants. In
Appendix \ref{sec:InvGath}  we will list all of the necessary invariants
which are obtained from Gathmann's program.

To obtain more relations on $a_{2}, \cdots a_{30}$ in equation \eqref{eqn:g3TRR}, we will
compute derivatives of $\Phi(\tau_{m, \alpha})$ at the origin of the big phase space
$t=0$.
From degree 0 part of $\Phi(\tau_{1, 1}) \mid _{t=0} = 0$, we obtain
\begin{equation} 0=  \frac{31}{967680} - \frac{a_{13}}{13824} - \frac{a_{14}}{13824}.
\end{equation}
From degree 0 part of $\Phi(\tau_{2, 0}) \mid _{t=0} = 0$, we obtain
\begin{equation} 0 =  -\frac{41}{290304} + \frac{7 a_{2}}{960} + \frac{a_{7}}{288} +
      \frac{a_{10}}{288} + \frac{a_{13}}{6912} + \frac{a_{16}}{288}.
      \end{equation}
From degree 1 part of $\Phi(\tau_{3, 1}) \mid _{t=0} = 0$, we obtain
\begin{eqnarray} 0 &=&  -\frac{1}{322560} - \frac{7 a_{2}}{2880} - \frac{a_{7}}{288} +
      \frac{a_{10}}{288} + \frac{a_{13}}{13824} \nonumber \\
      && - \frac{a_{14}}{13824} - \frac{a_{16}}{288}
       - \frac{a_{20}}{6} +
      8 a_{23} + \frac{a_{25}}{6} + \frac{a_{29}}{12}.
      \end{eqnarray}
From degree 0 part of $\tau_{2,1} \Phi(\tau_{0, 0}) \mid _{t=0} = 0$, we obtain
\begin{equation}0 =  \frac{31}{96768} - \frac{7 a_{4}}{46080} - \frac{a_{12}}{13824}
    - \frac{a_{13}}{13824} -
      \frac{a_{14}}{13824}.
      \end{equation}
From degree 0 part of $\tau_{3,0} \Phi(\tau_{0, 0}) \mid _{t=0} = 0$, we obtain
\begin{eqnarray}0  &=&  -\frac{1501}{725760} + \frac{7 a_{2}}{720} +
      \frac{a_{4}}{1920} + \frac{7 a_{5}}{2880} + \frac{a_{7}}{288} + \frac{a_{10}}{288}
      + \frac{a_{12}}{6912}  \nonumber \\
      &&     + \frac{a_{13}}{6912} + \frac{a_{15}}{288} + \frac{a_{16}}{288} + \frac{a_{24}}{288}
       + \frac{7 a_{28}}{720}.
       \end{eqnarray}
From degree 0 part of $\tau_{2,0} \Phi(\tau_{1, 0}) \mid _{t=0} = 0$, we obtain
\begin{eqnarray}0  &=&  -\frac{2329}{1451520} + \frac{7 a_{2}}{240} + \frac{7 a_{3}}{960} +
      \frac{a_{6}}{288} + \frac{a_{7}}{96} + \frac{a_{8}}{288} + \frac{a_{9}}{288} + \frac{a_{10}}{96} +
      \frac{a_{13}}{1728} \nonumber \\
      && + \frac{a_{14}}{2304} + \frac{a_{16}}{96} + \frac{a_{17}}{288}
      + \frac{a_{18}}{288}.
      \end{eqnarray}
From degree 1 part of $\tau_{3,1} \Phi(\tau_{0, 1}) \mid _{t=0} = 0$, we obtain
\begin{eqnarray}0  &=&  -\frac{31}{96768} + \frac{7 a_{4}}{138240} + \frac{7 a_{5}}{2880} -
      \frac{a_{7}}{288} + \frac{a_{10}}{288} + \frac{a_{12}}{13824} + \frac{a_{13}}{13824}
      - \frac{a_{14}}{13824} \nonumber \\
      && -  \frac{a_{16}}{288} - \frac{a_{20}}{6} + 8 a_{23} - \frac{a_{24}}{288} + \frac{a_{25}}{6} +
      \frac{a_{27}}{12} + \frac{7 a_{28}}{1440} + \frac{a_{29}}{8} - \frac{a_{30}}{24}.
      \end{eqnarray}
From degree 1 part of $\tau_{4,0} \Phi(\tau_{0, 1}) \mid _{t=0} = 0$, we obtain
\begin{eqnarray}0  &=&
    \frac{277}{207360} + \frac{a_{2}}{360} - \frac{13 a_{4}}{46080} - \frac{a_{5}}{120} -
    \frac{a_{7}}{144} +
     \frac{a_{12}}{13824} + \frac{a_{13}}{13824} + \frac{a_{14}}{13824} - \frac{a_{15}}{144}
     \nonumber \\
      && -  \frac{a_{16}}{144} -
      \frac{a_{20}}{6} + 16 a_{23} - \frac{a_{24}}{144} + \frac{a_{25}}{2} +
      \frac{a_{27}}{6} - \frac{11 a_{28}}{720} + \frac{a_{29}}{3} + \frac{a_{30}}{12}.
      \end{eqnarray}
From degree 1 part of $\tau_{3,1} \Phi(\tau_{1, 0}) \mid _{t=0} = 0$, we obtain
\begin{eqnarray}0  &=&
    \frac{41}{193536} - \frac{7 a_{2}}{720} - \frac{7 a_{3}}{2880} - \frac{a_{6}}{288} -
      \frac{a_{7}}{96} - \frac{a_{8}}{288} + \frac{a_{9}}{288} + \frac{a_{10}}{96}
      + \frac{a_{13}}{4608} +
      \frac{a_{14}}{13824} \nonumber \\
      && - \frac{a_{16}}{96} - \frac{a_{17}}{288} - \frac{a_{18}}{288} -
      \frac{a_{19}}{6} - \frac{2 a_{20}}{3} + 40 a_{23} + \frac{2 a_{25}}{3} + \frac{a_{26}}{6} +
      \frac{a_{29}}{3} + \frac{a_{30}}{12}.
      \end{eqnarray}
From degree 1 part of $\tau_{2,1} \Phi(\tau_{1, 1}) \mid _{t=0} = 0$, we obtain
\begin{eqnarray}0 &=&  -\frac{1}{46080} + \frac{7 a_{3}}{960} - \frac{7 a_{4}}{46080} -
      \frac{a_{6}}{288} - \frac{a_{7}}{288} + \frac{a_{8}}{288} + \frac{a_{9}}{288} +
      \frac{a_{10}}{96} -
      \frac{a_{12}}{13824} \nonumber \\
      && + \frac{a_{13}}{6912} - \frac{a_{16}}{288} - \frac{a_{19}}{6}
      - \frac{a_{20}}{2} -
      \frac{a_{22}}{12} + 24 a_{23} + \frac{a_{25}}{6} + \frac{a_{29}}{8} + \frac{a_{30}}{24}.
      \end{eqnarray}
From degree 1 part of $\tau_{3,0} \Phi(\tau_{1, 1}) \mid _{t=0} = 0$, we obtain
\begin{eqnarray} 0  &=&  -\frac{83}{193536} - \frac{7 a_{2}}{720} - \frac{a_{3}}{40} +
      \frac{a_{4}}{1920} + \frac{7 a_{5}}{2880} - \frac{a_{6}}{144} - \frac{5 a_{7}}{288} -
      \frac{a_{8}}{144} + \frac{a_{10}}{288} +
      \frac{a_{12}}{6912} \nonumber \\
      && + \frac{7 a_{13}}{13824} + \frac{5 a_{14}}{13824} +
      \frac{a_{15}}{288} - \frac{5 a_{16}}{288} - \frac{a_{17}}{144} - \frac{a_{18}}{144} -
      \frac{a_{19}}{6} - \frac{2 a_{20}}{3} + \frac{a_{22}}{6} + 56 a_{23} \nonumber \\
      && + \frac{a_{24}}{288} +
      a_{25} - \frac{a_{26}}{6} + \frac{7 a_{28}}{720} + \frac{7 a_{29}}{12}.
      \end{eqnarray}
From degree 1 part of $\tau_{2,1} \Phi(\tau_{2, 0}) \mid _{t=0} = 0$, we obtain
\begin{eqnarray}0  &=&  -\frac{19}{138240} + \frac{a_{2}}{320} + \frac{7 a_{4}}{46080} +
      \frac{a_{6}}{144} - \frac{5 a_{7}}{288} - \frac{a_{9}}{144} + \frac{a_{10}}{288}
       + \frac{a_{12}}{13824} \nonumber \\
      && + \frac{a_{13}}{6912} - \frac{5 a_{16}}{288} + \frac{a_{19}}{3} - a_{20} - \frac{a_{21}}{6} +
      72 a_{23} + a_{25} + \frac{a_{29}}{2}.
      \end{eqnarray}
From degree 1 part of $\tau_{1,1} \Phi(\tau_{2, 1}) \mid _{t=0} = 0$, we obtain
\begin{eqnarray}0  &= & \frac{7 a_{2}}{960} - \frac{a_{7}}{288} + \frac{5 a_{10}}{288} +
      \frac{a_{13}}{3456} - \frac{a_{14}}{2304} - \frac{a_{16}}{288} \nonumber \\
      && - \frac{2 a_{20}}{3}
      - \frac{a_{21}}{12} +
      24 a_{23} + \frac{a_{25}}{6} + \frac{a_{29}}{6}.
      \end{eqnarray}
From degree 1 part of $\tau_{2,0} \Phi(\tau_{2, 1}) \mid _{t=0} = 0$, we obtain
\begin{eqnarray}0 &=&  -\frac{1}{15360} -
      \frac{a_{2}}{240} + \frac{7 a_{3}}{960} - \frac{7 a_{4}}{46080} - \frac{a_{6}}{288} -
      \frac{a_{7}}{96} + \frac{a_{8}}{288} - \frac{a_{9}}{288} - \frac{a_{10}}{288} -
      \frac{a_{12}}{13824} \nonumber \\
      && +  \frac{a_{14}}{2304} - \frac{5 a_{16}}{288} + \frac{a_{17}}{288} -
      \frac{a_{18}}{288} - \frac{2 a_{20}}{3} - \frac{a_{22}}{6} + 72 a_{23} + a_{25}
      + \frac{a_{29}}{2}.
      \end{eqnarray}
From degree 1 part of $\tau_{1,1} \Phi(\tau_{3, 0}) \mid _{t=0} = 0$, we obtain
\begin{eqnarray} 0  &=&  -\frac{1}{46080} - \frac{107 a_{2}}{2880} - \frac{a_{7}}{32} +
      \frac{a_{10}}{288} + \frac{7 a_{13}}{13824} + \frac{5 a_{14}}{13824} -
      \frac{a_{16}}{32} \nonumber \\
      && - \frac{5 a_{20}}{6} + \frac{a_{21}}{6} +
      56 a_{23} + \frac{5 a_{25}}{6} + \frac{7 a_{29}}{12}.
      \end{eqnarray}
From degree 1 part of $\tau_{1, 1} \tau_{2,1} \Phi(\tau_{1, 0}) \mid _{t=0} = 0$, we obtain
\begin{eqnarray} 0 &=&  -\frac{1}{9216} + \frac{7 a_{2}}{240} + \frac{7 a_{3}}{480}
    - \frac{7 a_{4}}{46080} - \frac{a_{6}}{144} - \frac{a_{7}}{72} + \frac{a_{9}}{48} +
    \frac{a_{10}}{16} - \frac{a_{12}}{13824} \nonumber \\
      && +  \frac{a_{13}}{1152} - \frac{a_{16}}{72} - \frac{a_{17}}{288} -
      \frac{a_{18}}{288} - \frac{5 a_{19}}{6} - \frac{19 a_{20}}{6} - \frac{a_{21}}{4} - \frac{a_{22}}{4}
       +  144 a_{23} \nonumber \\
      && + \frac{5 a_{25}}{6} +
      \frac{a_{26}}{6} + \frac{19 a_{29}}{24} + \frac{5 a_{30}}{24}.
      \end{eqnarray}
From degree 1 part of $\tau_{2, 0} \tau_{2,0} \Phi(\tau_{1, 1}) \mid _{t=0} = 0$, we obtain
\begin{eqnarray} 0 &=& -\frac{661}{161280} -
      \frac{a_{2}}{96} - \frac{179 a_{3}}{1440} + \frac{263 a_{4}}{69120}
      + \frac{7 a_{5}}{480} - \frac{7 a_{6}}{144} - \frac{7 a_{7}}{144} - \frac{5 a_{8}}{144}
      - \frac{5 a_{9}}{144} \nonumber \\
      && - \frac{5 a_{10}}{144} + \frac{a_{11}}{144} + \frac{a_{12}}{768}
       - \frac{a_{13}}{1152} +
      \frac{a_{15}}{48} - \frac{13 a_{16}}{144} - \frac{5 a_{17}}{144} - \frac{7 a_{18}}{144}
      - a_{19} - \frac{5 a_{20}}{3} \nonumber \\
      && - \frac{a_{22}}{3} + 528 a_{23} +
      \frac{a_{24}}{48} + \frac{19 a_{25}}{3} -
      a_{26} + \frac{7 a_{28}}{96} + \frac{7 a_{29}}{2} - \frac{a_{30}}{6}.
      \label{eqn:R29}
      \end{eqnarray}

{\bf Remark}:
We have also checked more than 18000 other combinations of derivatives of
$\Phi$ for the ${\mathbb CP}^{1}$ case using a Mathematica program, but the
relations obtained are just linear combinations of the relations
\eqref{eqn:R1} to \eqref{eqn:R29}. This explicitly verifies that these
relations are indeed consistent.

\subsection{Proof of Theorem~\ref{thm:g3TRR}}
It is straightforward to solve $a_{2}, \ldots, a_{30}$ from equations \eqref{eqn:R1} to
\eqref{eqn:R29}. The answers are
{\large
\begin{equation}  \begin{array}{lclclcl}
 a_{2} = - \frac{1}{252},    & \hspace{10pt} &
 a_{3} = \frac{13}{168}, & \hspace{10pt} &
 a_{4} = \frac{41}{21},  & \hspace{10pt} &
 a_{5} = - \frac{13}{168},  \\
 a_{6} = \frac{1}{280},      &&
 a_{7} = - \frac{23}{5040},  &&
 a_{8} = - \frac{47}{5040},   &&
 a_{9} = - \frac{5}{1008},  \\
 a_{10} = \frac{23}{504},   &&
 a_{11} = \frac{11}{140},  &&
 a_{12} = - \frac{4}{35},   &&
 a_{13} = \frac{2}{105},  \\
 a_{14} = \frac{89}{210},   &&
 a_{15} = - \frac{1}{210},  &&
 a_{16} = \frac{1}{140},   &&
 a_{17} = \frac{23}{140},  \\
 a_{18} = - \frac{3}{140},   &&
 a_{19} = - \frac{1}{4480},  &&
 a_{20} = \frac{13}{8064},   &&
 a_{21} = - \frac{1}{2240},  \\
 a_{22} =  \frac{41}{6720},   &&
 a_{23} = \frac{1}{53760}, &&
 a_{24} = - \frac{1}{210},  &&
 a_{25} = - \frac{1}{5760},  \\
 a_{26} = - \frac{1}{2688},   &&
 a_{27} = - \frac{1}{5040}, &&
 a_{28} = \frac{5}{42},   &&
 a_{29} = \frac{1}{3780},  \\
 a_{30} = \frac{1}{252}.  &&&&
 \end{array}   \label{eqn:Coeffg3TRR}
\end{equation}
}
Together with equations \eqref{eqn:a1} and \eqref{eqn:g3TRRa}, this proves
equation \eqref{eqn:g3TRR}.
$\Box$

\vspace{10pt}
We make the following observations about topological recursion relations of
genus less than or equal to 3:
\begin{enumerate}
\item $T(\vw)$ does not appear in the lower genus terms of the genus-$g$
topological recursion relations for $T^{g}(\vw)$ when $g=1, 2, 3$.
\item With the exception of $a_{25}$, all of the denominators in the coefficients of
  the genus-3 topological recursion relation \eqref{eqn:g3TRR} have a factor $7$.
We also note that the genus-2 topological recursion relations (i.e. Mumford's equation),
have denominators in its coefficients
of lower genus terms contain a factor of 10, while for the genus-1
  topological recursion relation, the corresponding factor is 24.
\item For $g=1, 2, 3$, the coefficients of the terms consisting of purely
  genus-0 functions in the genus-$g$ topological recursion relations are
\[ \frac{1}{(2g+1)!! \,\, 8^{g}}.\]
We conjecture that this should also hold for all genera.
\end{enumerate}

\section{A new relation in tautological ring of $\overline{\cal M}_{3,1}$}
\label{sec:taut}

Note that equation \eqref{eqn:g3TRRa} is a direct translation of a relation
in the tautological ring of $\overline{\cal M}_{3, 1}$ with undetermined
coefficients $a_{1}, \cdots a_{30}$. Those coefficients were determined
during the proof of topological recursion relation \eqref{eqn:g3TRR}. Therefore we
have also obtained a proof to the following theorem:

\begin{thm}
In the tautological ring of $\overline{\cal M}_{3, 1}$, the following relation holds
{\allowdisplaybreaks
\begin{eqnarray}
\psi_{1}^{3} &=&  - \frac{1}{126} \hspace{10pt}
\begin{picture}(80, 30)
\put(7, 5){\line(0, 1){10}}
\put(5.5, -1){$\scriptscriptstyle 2$}
\put(7, 1){\circle{7}}
\put(40.5, 2){\vector(-1, 0){30}}
\put(40, 1){ \circle{7}}
\curve(48, 2, 60, 9, 70, 2, 60, -7, 48, 0)
\end{picture}
+ \frac{13}{84} \hspace{10pt}
\begin{picture}(80, 30)
\put(5.5, -1){$\scriptscriptstyle 2$}
\put(7, 1){\circle{7}}
\put(40.5, 2){\vector(-1, 0){30}}
\put(40, 1){ \circle{7}}
\curve(48, 2, 60, 9, 70, 2, 60, -7, 48, 0)
\put(43, 5){\line(0, 1){10}}
\end{picture}
+\frac{5}{42} \hspace{10pt}
\begin{picture}(80, 30)
\put(5.5, -1){$\scriptscriptstyle 2$}
\put(7, 1){\circle{7}}
\put(40.5, 2){\vector(-1, 0){30}}
\put(40, 1){ \circle{7}}
\curve(11, -1, 26, -6, 41, -1)
\put(43, 5){\line(0, 1){10}}
\end{picture} \nonumber \\
&& +\frac{41}{21} \hspace{10pt}
\begin{picture}(100, 30)
\put(5.5, -1){$\scriptscriptstyle 2$}
\put(7, 1){\circle{7}}
\put(40.5, 2){\vector(-1, 0){30}}
\put(40, 1){ \circle{7}}
\put(43, 5){\line(0, 1){10}}
\put(48, 1){\line(1, 0){30}}
\put(81, 1){\circle{7}}
\put(79.5, -1){$\scriptscriptstyle 1$}
\end{picture}
- \frac{13}{84} \hspace{10pt}
\begin{picture}(120, 30)
\put(5.5, -1){$\scriptscriptstyle 2$}
\put(7, 1){\circle{7}}
\put(10.5, 2){\line(1, 0){30}}
\put(40, 1){ \circle{7}}
\put(43, 5){\line(0, 1){10}}
\put(48, 1){\line(1, 0){30}}
\put(81, 1){\circle{7}}
\curve(84.5, 2, 97, 9, 107, 2, 97, -7, 84.5, 0)
\end{picture} \nonumber \\
&& + \frac{1}{140} \hspace{10pt}
\begin{picture}(120, 30)
\put(5.5, -1){$\scriptscriptstyle 1$}
\put(7, 1){\circle{7}}
\put(7, 5){\line(0, 1){10}}
\put(10.5, 2){\line(1, 0){30}}
\put(40, 1){ \circle{7}}
\put(42.5, -1){$\scriptscriptstyle 1$}
\put(48, 1){\line(1, 0){30}}
\put(81, 1){\circle{7}}
\curve(85, 2, 97, 9, 107, 2, 97, -7, 85, 0)
\end{picture}
- \frac{23}{2520} \hspace{10pt}
\begin{picture}(120, 30)
\put(5.5, -1){$\scriptscriptstyle 1$}
\put(7, 1){\circle{7}}
\put(10.5, 2){\line(1, 0){30}}
\put(40, 1){ \circle{7}}
\put(42.5, -1){$\scriptscriptstyle 1$}
\put(43, 5){\line(0, 1){10}}
\put(48, 1){\line(1, 0){30}}
\put(81, 1){\circle{7}}
\curve(85, 2, 97, 9, 107, 2, 97, -7, 85, 0)
\end{picture}  \nonumber \\
&& - \frac{47}{2520} \hspace{10pt}
\begin{picture}(120, 30)
\put(5.5, -1){$\scriptscriptstyle 1$}
\put(7, 1){\circle{7}}
\put(10.5, 2){\line(1, 0){30}}
\put(40, 1){ \circle{7}}
\put(42.5, -1){$\scriptscriptstyle 1$}
\put(48, 1){\line(1, 0){30}}
\put(81, 1){\circle{7}}
\curve(85, 2, 97, 9, 107, 2, 97, -7, 85, 0)
\put(81, 5){\line(0, 1){10}}
\end{picture}
- \frac{1}{105} \hspace{10pt}
\begin{picture}(100, 30)
\put(5.5, -1){$\scriptscriptstyle 1$}
\put(7, 1){\circle{7}}
\put(10.5, 2){\line(1, 0){30}}
\put(43.5, -1){$\scriptscriptstyle 1$}
\put(44, 1){\circle{7}}
\curve(48, 2, 66, 6, 81, 2)
\curve(48, -2, 66, -6, 81, -2)
\put(80, 1){ \circle{7}}
\put(82, 5){\line(0, 1){10}}
\end{picture} \nonumber \\
&& +\frac{11}{70} \hspace{10pt}
\begin{picture}(100, 30)
\put(5.5, -1){$\scriptscriptstyle 1$}
\put(7, 1){\circle{7}}
\put(10.5, 2){\line(1, 0){30}}
\put(40, 1){ \circle{7}}
\put(43, 5){\line(0, 1){10}}
\put(48, 1){\line(1, 0){30}}
\put(81, 1){\circle{7}}
\put(79.5, -1){$\scriptscriptstyle 1$}
\curve(7, -2, 22, -8, 44, -10, 66, -8, 81, -2)
\end{picture}
- \frac{4}{35} \hspace{10pt}
\begin{picture}(100, 30)
\put(5.5, -1){$\scriptscriptstyle 1$}
\put(7, 1){\circle{7}}
\put(10.5, 2){\line(1, 0){30}}
\put(42.5, -1){$\scriptscriptstyle 1$}
\put(44, 1){\circle{7}}
\put(47, 2){\line(1, 0){30}}
\put(76, 1){ \circle{7}}
\put(79, 5){\line(0, 1){10}}
\put(84, 1){\line(1, 0){30}}
\put(118, 1){\circle{7}}
\put(116.5, -1){$\scriptscriptstyle 1$}
\end{picture} \nonumber \\
&& - \frac{1}{105} \hspace{10pt}
\begin{picture}(120, 30)
\put(5.5, -1){$\scriptscriptstyle 1$}
\put(7, 1){\circle{7}}
\put(10.5, 2){\line(1, 0){30}}
\put(40, 1){ \circle{7}}
\put(43, 5){\line(0, 1){10}}
\put(48, 1){\line(1, 0){30}}
\put(81, 1){\circle{7}}
\put(80.5, -1){$\scriptscriptstyle 1$}
\curve(85, 2, 97, 9, 107, 2, 97, -7, 85, 0)
\end{picture}
+ \frac{1}{70} \hspace{10pt}
\begin{picture}(100, 30)
\put(5.5, 5){\line(0, 1){10}}
\put(5.5, -1){$\scriptscriptstyle 1$}
\put(7, 1){\circle{7}}
\curve(11, 2, 27, 6, 44, 2)
\curve(11, -2, 27, -6, 44, -2)
\put(43, 1){ \circle{7}}
\put(50.5, 2){\line(1, 0){30}}
\put(83.5, -1){$\scriptscriptstyle 1$}
\put(85, 1){\circle{7}}
\end{picture} \nonumber \\
&& + \frac{23}{70} \hspace{10pt}
\begin{picture}(100, 30)
\put(5.5, -1){$\scriptscriptstyle 1$}
\put(7, 1){\circle{7}}
\curve(11, 2, 27, 6, 44, 2)
\curve(11, -2, 27, -6, 44, -2)
\put(43, 1){ \circle{7}}
\put(46, 5){\line(0, 1){10}}
\put(50.5, 2){\line(1, 0){30}}
\put(83.5, -1){$\scriptscriptstyle 1$}
\put(85, 1){\circle{7}}
\end{picture}
- \frac{3}{70} \hspace{10pt}
\begin{picture}(100, 30)
\put(5.5, -1){$\scriptscriptstyle 1$}
\put(7, 1){\circle{7}}
\curve(11, 2, 27, 6, 44, 2)
\curve(11, -2, 27, -6, 44, -2)
\put(43, 1){ \circle{7}}
\put(50.5, 2){\line(1, 0){30}}
\put(83.5, -1){$\scriptscriptstyle 1$}
\put(85, 1){\circle{7}}
\put(85, 5){\line(0, 1){10}}
\end{picture}
\nonumber \\
&& - \frac{5}{504} \hspace{10pt}
\begin{picture}(100, 30)
\put(7, 5){\line(0, 1){10}}
\put(5.5, -1){$\scriptscriptstyle 1$}
\put(7, 1){\circle{7}}
\put(10.5, 2){\line(1, 0){30}}
\put(40, 1){ \circle{7}}
\curve(41, -2, 38, -8, 44, -16, 50, -8, 47, -2)
\put(48, 1){\line(1, 0){30}}
\put(81, 1){\circle{7}}
\put(79.5, -1){$\scriptscriptstyle 1$}
\end{picture}
+ \frac{23}{126} \hspace{10pt}
\begin{picture}(100, 30)
\put(5.5, -1){$\scriptscriptstyle 1$}
\put(7, 1){\circle{7}}
\put(10.5, 2){\line(1, 0){30}}
\put(40, 1){ \circle{7}}
\curve(41, -2, 38, -8, 44, -16, 50, -8, 47, -2)
\put(43, 5){\line(0, 1){10}}
\put(48, 1){\line(1, 0){30}}
\put(81, 1){\circle{7}}
\put(79.5, -1){$\scriptscriptstyle 1$}
\end{picture} \nonumber \\
&& + \frac{4}{105} \hspace{10pt}
\begin{picture}(100, 30)
\put(7, 5){\line(0, 1){10}}
\put(5.5, -1){$\scriptscriptstyle 1$}
\put(7, 1){\circle{7}}
\put(10.5, 2){\line(1, 0){30}}
\put(40, 1){ \circle{7}}
\put(43.5, -2.5){\line(0, -1){20}}
\put(43.5, -26.5){\circle{7}}
\put(42.5, -28.5){$\scriptscriptstyle 1$}
\put(48, 1){\line(1, 0){30}}
\put(81, 1){\circle{7}}
\put(79.5, -1){$\scriptscriptstyle 1$}
\end{picture}
+ \frac{89}{35} \hspace{10pt}
\begin{picture}(100, 30)
\put(5.5, -1){$\scriptscriptstyle 1$}
\put(7, 1){\circle{7}}
\put(10.5, 2){\line(1, 0){30}}
\put(43, 5){\line(0, 1){10}}
\put(40, 1){ \circle{7}}
\put(43.5, -2.5){\line(0, -1){20}}
\put(43.5, -26.5){\circle{7}}
\put(42.5, -28.5){$\scriptscriptstyle 1$}
\put(48, 1){\line(1, 0){30}}
\put(81, 1){\circle{7}}
\put(79.5, -1){$\scriptscriptstyle 1$}
\end{picture} \nonumber \\
&& + \frac{1}{630} \hspace{10pt}
\begin{picture}(60, 50)
\put(5.5, -1){$\scriptscriptstyle 1$}
\put(7, 1){\circle{7}}
\put(7, 5){\line(0, 1){10}}
\curve(11, 2, 27, 6, 44, 2)
\put(10, 0){\line(1, 0){32}}
\curve(11, -2, 27, -6, 44, -2)
\put(43, 1){ \circle{7}}
\end{picture}
+ \frac{1}{42} \hspace{10pt}
\begin{picture}(60, 50)
\put(5.5, -1){$\scriptscriptstyle 1$}
\put(7, 1){\circle{7}}
\curve(11, 2, 27, 6, 44, 2)
\put(10, 0){\line(1, 0){32}}
\curve(11, -2, 27, -6, 44, -2)
\put(43, 1){ \circle{7}}
\put(46, 5){\line(0, 1){10}}
\end{picture}
- \frac{1}{1260} \hspace{10pt}
\begin{picture}(100, 50)
\curve(28.5, 2, 15, 8, 2, 0, 15, -8, 30, -2)
\put(31, -1){$\scriptscriptstyle 1$}
\put(32.5, 1){\circle{7}}
\curve(37, 2, 56, 6, 71, 2)
\curve(37, -2, 56, -6, 71, -2)
\put(70, 1){ \circle{7}}
\put(73, 5){\line(0, 1){10}}
\end{picture} \nonumber \\
&& - \frac{1}{560} \hspace{10pt}
\begin{picture}(100, 30)
\curve(28.5, 2, 15, 8, 2, 0, 15, -8, 30, -2)
\put(32.5, 1){\circle{7}}
\curve(36, 2, 56, 6, 71, 2)
\curve(35, -2, 56, -6, 71, -2)
\put(70, 1){ \circle{7}}
\put(72.5, -1){$\scriptscriptstyle 1$}
\put(73, 5){\line(0, 1){10}}
\end{picture}
+ \frac{41}{1680} \hspace{10pt}
\begin{picture}(100, 30)
\curve(28.5, 2, 15, 8, 2, 0, 15, -8, 30, -2)
\put(32.5, 1){\circle{7}}
\put(32, 5){\line(0, 1){10}}
\curve(36, 2, 56, 6, 71, 2)
\curve(35, -2, 56, -6, 71, -2)
\put(70, 1){ \circle{7}}
\put(72.5, -1){$\scriptscriptstyle 1$}
\end{picture} \nonumber \\
&& - \frac{1}{1440} \hspace{10pt}
\begin{picture}(110, 30)
\curve(28.5, 2, 15, 8, 2, 0, 15, -8, 30, -2)
\put(32.5, 1){\circle{7}}
\put(31, -1){$\scriptscriptstyle 1$}
\put(36, 2){\line(1, 0){30}}
\put(32, 5){\line(0, 1){10}}
\put(66, 1){ \circle{7}}
\curve(73, 2, 87, 9, 98, 2, 87, -7, 73, 0)
\end{picture}
- \frac{1}{672} \hspace{10pt}
\begin{picture}(110, 30)
\curve(28.5, 2, 15, 8, 2, 0, 15, -8, 30, -2)
\put(32.5, 1){\circle{7}}
\put(31, -1){$\scriptscriptstyle 1$}
\put(36, 2){\line(1, 0){30}}
\put(69, 5){\line(0, 1){10}}
\put(66, 1){ \circle{7}}
\curve(73, 2, 87, 9, 98, 2, 87, -7, 73, 0)
\end{picture} \nonumber \\
&& - \frac{1}{560} \hspace{10pt}
\begin{picture}(80, 30)
\put(7, 5){\line(0, 1){10}}
\put(5.5, -1){$\scriptscriptstyle 1$}
\put(7, 1){\circle{7}}
\put(10.5, 2){\line(1, 0){30}}
\put(40, 1){ \circle{7}}
\curve(44.5, 5, 52, 14, 64, 18, 54, 6, 47, 3)
\curve(43, -3, 52, -14, 64, -18, 54, -6, 46, -2)
\end{picture}
+ \frac{13}{1008} \hspace{10pt}
\begin{picture}(80, 30)
\put(5.5, -1){$\scriptscriptstyle 1$}
\put(7, 1){\circle{7}}
\put(10.5, 2){\line(1, 0){30}}
\put(44, 1){\circle{7}}
\put(44, 5){\line(0, 1){10}}
\curve(44.5, 5, 52, 14, 64, 18, 54, 6, 47, 3)
\curve(43, -3, 52, -14, 64, -18, 54, -6, 46, -2)
\end{picture}
 + \frac{1}{1120} \hspace{10pt}
\begin{picture}(80, 30)(24, 0)
\put(44, 1){\circle{7}}
\put(41.5, 3){\line(-1, 1){10}}
\curve(44.5, 5, 52, 14, 64, 18, 54, 6, 47, 3)
\curve(43.5, -3, 36, -12, 24, -16, 34, -4, 41, -1)
\curve(43, -3, 52, -14, 64, -18, 54, -6, 46, -2)
\end{picture} \label{eqn:Graphg3TRR}
\end{eqnarray}}
\end{thm}
In this formula, each stratum in $\overline{\cal M}_{g, n}$ is represented by its dual
graph. We adopt the conventions of \cite{Ge2} for dual graphs with a slight
modification. We denote vertices of genus $0$  by a hollow circle
$\begin{picture}(10,8)(0,0)\put(5,3){\circle{7}}\end{picture}$, and
vertices of genus $g \geq 1$ by $\begin{picture}(10,8)(0,0)\put(5,3){\circle{7}}
\put(3,2){$\scriptscriptstyle g$}\end{picture}$.
A vertex with an incident arrowhead denotes the $\psi$ class associated to
the marked point (which is at a node) on the irreducible component associated
to that vertex.

Note that when translating relations in the tautological ring of $\overline{\cal M}_{g, n}$
to Gromov-Witten invariants, we need to divide the coefficient of each stratum by the number of
elements in the automorphism group of the corresponding dual graph. This
explains the difference between the coefficients in this formula and those in equation
\eqref{eqn:g3TRR}.

\section{Application to Higher Spin Curves}

We briefly review the moduli space of $r$-spin curves and $r$-spin
theory. For details, we refer the reader to \cite{JKV-99}. Let
$r\geq 2$ be an integer.  For each $r$, $r$-spin theory is a
cohomological field theory, in the sense of Kontsevich-Manin, just
as is the Gromov-Witten theory of a compact, symplectic manifold
$M$. Consequently, the correlation functions in this theory satisfy
the same universal equations as does Gromov-Witten theory. In
particular, they satisfy the topological recursion relations, the
string, and dilaton equations. However, there is no analog of the
divisor equation but we will not need this.

The role of the moduli space of stable maps is replaced by the
moduli space of $r$-spin curves $\overline{\cal M}_{g,n}^{1/r}$. Let
$[n] := \{\,0,\ldots,n \}$ for any nonnegative integer $n$ then we
have the disjoint union $\overline{\cal M}_{g,n}^{1/r} :=
\bigsqcup_{\ba\in [r-2]^n} \overline{\cal M}_{g,n}^{1/r}(\ba)$.
While the definition of an $r$-spin curve is rather involved, over a
smooth stable curve $(C;p_1,\ldots,p_n)$ of genus $g$ with $n$
marked points when $r$ is prime, a point in $\overline{\cal
M}_{g,n}^{1/r}(\ba)$ consists of a line bundle $\cl\to C$ together
with an isomorphism $\cl^{\otimes r}\to \omega(-\sum_{i=1}^n \alpha_i p_i)$
where $\ba = (\alpha_1,\ldots,\alpha_n)$ and $\omega$ is the
canonical bundle of $C$. For degree reasons, such a line bundle
exists only if $(2-2g+\sum_{i=1}^n \alpha_i)/r$ is an integer. When
$r$ is not prime then one should introduce all $d$-th roots for all
$d$ which divides $r$. This moduli space admits a compactification
which is a smooth Deligne-Mumford stack by allowing the bundles to
degenerate over the nodes. There are also the tautological classes
$\Psi_i$ associated to the cotangent line at the marked points
$p_i$. When nonempty, the dimension of  $\overline{\cal
M}_{g,n}^{1/r}(\ba)$ is the same as that of $\M_{g,n}$, namely
$3g-3+n$, since there are only a finite number of $r$-spin
structures over a given stable curve. There is no analog in $r$-spin
theory of the degree $d$ of a stable map into a manifold so we will
ignore this parameter.

The virtual fundamental class in Gromov-Witten theory is replaced by the
virtual class $c^{1/r}$, a cohomology class in $H^{2 D}(\overline{\cal
  M}^{1/r}_{g,n}(\ba))$ where $D = \frac{1}{r}((r-2)(g-1)+\sum_{i=1}^n \alpha_i)$.
The analog of the cohomology ring of $M$ is the vector space ${\cal H}$ with a
basis $\{\,\gamma_0\,\ldots\,\gamma_{r-2}\,\}$ with the pairing $\eta_{\alpha
  \beta} = \delta_{\alpha+\beta,r-2}$. The correlation functions are defined by
\[ \gwih{g}{\grav{n_{1}}{\alpha_{1}} \, \grav{n_{2}}{\alpha_{2}} \,
    \ldots \, \grav{n_{k}}{\alpha_{k}}} :=
    r^{1-g} \int_{[\overline{\cal M}^{1/r}_{g,k}({\boldsymbol \alpha})]}
    c^{1/r}\cup\bigcup_{i=1}^k \Psi^{n_i}
\]
where $[\overline{\cal M}^{1/r}_{g,k}({\boldsymbol \alpha})]$ denotes the
fundamental class of the stack $\overline{\cal M}^{1/r}_{g,k}({\boldsymbol
  \alpha})$.

As far as dimensional considerations are concerned, $r$-spin theory behaves
as if it were Gromov-Witten theory of a manifold with real dimension $2(r-2)/r$ and
$\gamma_\alpha$ behaves as if it were a cohomology class of dimension
$2\alpha/r$ for all $\alpha = 0,\dots,r-2$.

Applying Theorem \ref{thm:g3TRR} and dimensional arguments, we obtain the
following.

\begin{pro}
Let $r\geq 2$. In $r$-spin theory, the only possible nontrivial $1$-point
genus-$3$ correlators are as follows:

\begin{description}
\item[$r=2$:] $\left<\tau_{7,0}\right>_3 = \frac{1}{82944}$
\item[$r=3$:] $\left<\tau_{6,1}\right>_3 = \frac{1}{31104}$
\item[$r=4$:] $\left<\tau_{6,0}\right>_3 = \frac{3}{20480}$
\item[$r=5$:] All $1$-point genus $3$-correlators vanish for dimensional reasons.
\item[$r=6$:] $\left<\tau_{5,4}\right>_3 =\frac{2561}{20901888}$
\item[$r\geq 7$:] $\left<\tau_{5,4}\right>_3$. Can be reduced to 5-point
  genus-0 correlators which can be calculated using the WDVV equation.
\end{description}

Furthermore, when $r=3$, the only possible nontrivial $2$-point correlators are

\[
\begin{array}{lllll}
\gwih{3}{\tau_{7,0}\tau_{0,1}} = 1/15552, &\hspace{15pt} &
\gwih{3}{\tau_{6,0} \tau_{1,1}} = 19/77760, &\hspace{15pt} &
\gwih{3}{\tau_{5,0}\tau_{2,1}} = 47/77760,  \\
\gwih{3}{\tau_{4,0}\tau_{3,1}} = 67/77760,  &&
\gwih{3}{\tau_{3,0}\tau_{4,1}} = 443/77760, &&
\gwih{3}{\tau_{2,0}\tau_{5,1}} = 103/217728, \\
\gwih{3}{\tau_{1, 0} \tau_{6,1}} = 5/31104, &&
\gwih{3}{\tau_{0,0} \tau_{7,1}} = 1/31104. && \\
\end{array}
\]
\end{pro}

\begin{rem}
When $r=3$, our results for $\gwih{3}{\tau_{7,0}\tau_{0,1}}$ and
$\gwih{3}{\tau_{6,1}}$ both agree with the results of Shadrin
\cite{Sh-02} who calculated them by studying the geometry of
$\M^{1/3}_{3,n}$.  He also showed that both of these numbers were consistent
with the $3$-rd KdV hierarchy as predicted by the generalized Witten conjecture.
\end{rem}

\appendix
\vspace{30pt}
\hspace{160pt} {\bf \Large Appendix}

\section{Gromov-Witten invariants of ${\mathbb CP}^{1}$ used to
determine the genus-3 topological recursion relation}
\label{sec:InvGath}

We need the following Gromov-Witten invariants of ${\mathbb CP}^{1}$
which are obtained using Gathmann's program based on the Virasoro constraints.

Genus-1 invariants:
\allowdisplaybreaks
\[
\gwih{1,1}{\tau_{1, 1}^{2}} = \gwih{1,1}{\tau_{3, 0}}= 0, \hspace{15pt}
\gwih{1,1}{\tau_{1, 1} \tau_{2, 0}} = 1/8, \hspace{15pt}
\gwih{1,1}{\tau_{2, 0}^{2}} = -1/6, \hspace{15pt}
\gwih{1,1}{\tau_{2, 1}} = 1/24.
\]

Genus-2 invariants:
\[
\begin{array}{lllll}
\gwih{2,0}{\tau_{1, 1}^{2}} = 0, &\hspace{15pt} &
\gwih{2,0}{\tau_{1, 1} \tau_{2, 0}} = 7/1920, &\hspace{15pt} &
\gwih{2,0}{\tau_{2, 0}^{2}} = -5/288,  \\
\gwih{2,0}{\tau_{2, 1}} = 7/5760,  &&
\gwih{2,0}{\tau_{3, 0}} = -1/240, &&
\gwih{2,1}{\tau_{1, 1}^{4}} = 0, \\
\gwih{2,1}{\tau_{1, 1}^{3} \tau_{2, 0}} = 0, &&
\gwih{2,1}{\tau_{1, 1}^{2} \tau_{2, 0}^{2}} = 1/32, &&
\gwih{2,1}{\tau_{1, 1} \tau_{2, 0}^{3}} = 17/64, \\
\gwih{2,1}{\tau_{2, 0}^{4}} = 5/6, &&
\gwih{2,1}{\tau_{1, 1}^{2} \tau_{2, 1}} = 0, &&
\gwih{2,1}{\tau_{1, 1} \tau_{2, 0} \tau_{2, 1}} = 1/192, \\
\gwih{2,1}{\tau_{2, 0}^{2} \tau_{2, 1}} = 23/576, &&
\gwih{2,1}{\tau_{2, 1}^{2}} = 1/576, &&
\gwih{2,1}{\tau_{1, 1}^{2} \tau_{3, 0}} = 0, \\
\gwih{2,1}{\tau_{1, 1} \tau_{2, 0} \tau_{3, 0}} = 1/32, &&
\gwih{2,1}{\tau_{2, 0}^{2} \tau_{3, 0}} = 11/96, &&
\gwih{2,1}{\tau_{2, 1} \tau_{3, 0}} = 1/192, \\
\gwih{2,1}{\tau_{3, 0}^{2}} = 29/1440, &&
\gwih{2,1}{\tau_{1, 1} \tau_{3, 1}} = 0, &&
\gwih{2,1}{\tau_{2, 0} \tau_{3, 1}} = 1/192, \\
\gwih{2,1}{\tau_{1, 1} \tau_{4, 0}} = 1/384, &&
\gwih{2,1}{\tau_{2, 0} \tau_{4, 0}} = 41/2880, &&
\gwih{2,1}{\tau_{4, 1}} = 1/1920, \\
\gwih{2,1}{\tau_{5, 0}} = 1/576. &&&&
\end{array}
\]

Genus-3 invariants
\[
\begin{array}{lllll}
\gwih{3, 0}{\tau_{1, 1}^{2} \tau_{2, 1}}=0, & \hspace{1pt} &
\gwih{3,0}{\tau_{1, 1} \tau_{2, 0} \tau_{2, 1}}=0, & \hspace{1pt} &
\gwih{3, 0}{\tau_{2, 0}^{2} \tau_{2, 1}}=-31/10752, \\
\gwih{3, 0}{\tau_{2, 1}^{2}}=0, &&
\gwih{3,0}{\tau_{2, 1} \tau_{3, 0}}= -31/96768, &&
\gwih{3, 0}{\tau_{3, 0}^{2}}=1501/725760, \\
\gwih{3,0}{\tau_{1, 1} \tau_{3, 1}}=0, &&
\gwih{3,0}{\tau_{2, 0} \tau_{3, 1}}=-31/96768, &&
\gwih{3,0}{\tau_{1, 1} \tau_{4, 0}}= -31/193536, \\
\gwih{3,0}{\tau_{4, 1}}= -31/967680, &&
\gwih{3,0}{\tau_{5, 0}} = 41/290304, &&
\gwih{3, 1}{\tau_{3, 1}^{2}}=0, \\
\gwih{3,1}{\tau_{1, 1} \tau_{2, 1} \tau_{4, 0}}=1/9216, &&
\gwih{3,1}{\tau_{3, 1} \tau_{4, 0}} = 1/9216, &&
\gwih{3, 1}{\tau_{2, 0}^{2} \tau_{4, 1}} = 7/5760, \\
\gwih{3,1}{\tau_{2, 1} \tau_{4, 1}} = 1/46080, &&
\gwih{3,1}{\tau_{3, 0} \tau_{4, 1}} = 1/9216, &&
\gwih{3,1}{\tau_{2, 1} \tau_{5, 0}}=19/138240, \\
\gwih{3,1}{\tau_{1, 1} \tau_{5, 1}}=0, &&
\gwih{3,1}{\tau_{2, 0} \tau_{5, 1}}=1/15360, &&
\gwih{3,1}{\tau_{1, 1} \tau_{6, 0}}=1/46080, \\
\gwih{3,1}{\tau_{6, 1}}=1/322560. &&&&
\end{array}
\]

All other invariants needed to determine the genus-3 topological recursion
relation can be computed from the string, dilaton, divisor equations, and the
0-point invariants \[ \gwih{0, 1}{} = 1, \hspace{20pt} \gwih{1,0}{}=0.\]


\vspace{30pt} \noindent
Takashi Kimura \\
Department of Mathematics and Statistics \\
111 Cummington St. \\
Boston University \\
Boston, MA 02215, USA \\
E-mail address: {\it kimura@math.bu.edu}

\vspace{30pt} \noindent
Xiaobo Liu \\
Department of Mathematics  \\
University of Notre Dame \\
Notre Dame,  IN  46556, USA \\
E-mail address: {\it xliu3@nd.edu}

\end{document}